\documentclass[journal]{article}

\usepackage{graphicx,url}
\usepackage{dcolumn}
\usepackage{bm}
\usepackage{amsmath,amsfonts,amssymb}

\newtheorem{definition}{\noindent{\it Definition}}[section]
\newtheorem{theorem}{\noindent{\it Theorem}}[section]

\begin{document}

\title{On Finitely Stable Additive Bases}

\author{Lucas Y. Obata, Luan A. Ferreira, Giuliano G. La Guardia
\thanks{Lucas Obata and Giuliano La Guardia (corresponding author) are with the Department of Mathematics and Statistics,
State University of Ponta Grossa (UEPG), 84030-900, Ponta Grossa,
PR, Brazil, e-mails: (obatalucas@gmail.com;gguardia@uepg.br). Luan Ferreira is with Institute of Mathematics and
Statistics of University of S\~ao Paulo, Rua do Mat\~ao, 1010, Cidade Universit\'aria, S\~ao  Paulo, SP 05508-090, Brazil,
e-mail: (luan@ime.usp.br)}}

\maketitle

\begin{abstract}
The concept of additive basis has been investigated in the literature for several 
mathematicians which works with number theorem. Recently, the concept of 
finitely stable additive basis was introduced. In this note we provide a counterexample of the reciprocal of
Theorem 2.2 shown in [Ferreira, L.A.: Finite Stable Additive Basis; Bull. Aust. Math. Soc.].
The idea of the construction of such a counterexample can possibly help the process of
finding additive bases of some specific order.
\end{abstract}

\section{Introduction}\label{sec1}

The concept of additive basis as well as finitely stable (additive) basis are well-known
and much investigated for number theory theorists. For more details about such a theory,
we refer the reader to the textbook \cite{Nathanson:1996}; see also the interesting papers 
\cite{Nathanson:1972,Nathanson:1974,Nathanson:1996A,Nathanson:2005,Nathanson:2012,Nathanson:2014,Tafula:2019}.

An additive basis $A$ is a subset $A\subseteq {\mathbb N}=\{0, 1, 2, 3, \ldots \}$ having the property
that there exists a positive integer $h$ such that every $n \in {\mathbb N}$ can
be written as the sum of $h$ elements of $A$. The order $\operatorname{o}(A)$ of $A$ is the minimum
among all $h$ satisfying such condition.
From Lagrange's Theorem it follows
that the set of squares
$${\mathbb N}^2 =\{ 0, 1, 4, 9, \ldots, n^2, \ldots\}$$
is an additive basis of order $4$.
Applying Wieferich's Theorem it follows that the set of cubes
$${\mathbb N}^3 = \{0, 1, 8, \ldots, n^3 , \ldots \}$$
has order $9$.

Let us recall the concept of finitely stable basis.

\begin{definition}\label{fsbasis}
Let $A$ be an additive basis. Then $A$ is said to be finitely stable if
$\operatorname{o}(A)=\operatorname{o}(A\cup F)$ for every finite set 
$F \subset {\mathbb N}$.
\end{definition}

\section{The Counterexample}\label{sec2}

We maintain the notation adopted in \cite{Ferreira:2018}. Let $A, B \subseteq 
{\mathbb N}$. We define
$$ A + B = \{a + b : a \in A, b \in B\}.$$
If $t$ is a positive integer and $A \subseteq {\mathbb N}$ we write
$$tA := \underbrace{A + A + \cdots + A}_{ t \ \mbox{times}}$$ and $0A = \{0\}$.
If $n \in {\mathbb N}$, we denote $A(n)=| \{a \in A : 1 \leq a\leq n \} |$.

One of us provided in Ref.~\cite{Ferreira:2018} conditions in order to guarantee 
finitely stability of a given (additive) basis.

\begin{theorem}\cite[Theorem 2.2.]{Ferreira:2018}\label{LF}
Let $A$ be an additive basis whose order satisfies $\operatorname{o}(A)=h \geq 3$. If
$$\displaystyle\lim_{n\rightarrow \infty} \frac{((h-2)A)(n)}{n}=0$$ and
$$\displaystyle\limsup_{n\rightarrow \infty} \frac{((h-1)A)(n)}{n}< 1,$$ then $A$ is finitely stable.
\end{theorem}

We here present a counterexample for the reciprocal of Theorem~\ref{LF}.

\vspace{0.3cm}

\textbf{Counterexample}

\vspace{0.3cm}

Let

$A = \{0, 1, \ldots , 10\}\cup \{22, 23, \ldots, 100\}\cup \ldots \cup 
\{ 2\cdot 10^{n-1}+2, 2 \cdot 10^{n-1}+3, \ldots , 10^{n}\}\cup \ldots =
A_1 \cup A_2 \cup \ldots \cup A_n \cup \ldots$.

\vspace{0.5cm}

It is easy to see that $A$ is an additive basis of order 3, i.e., all nonnegative integers can be written
as sum of three elements of $A$. Moreover, $A$ is finitely stable since, for
every finite set $F \subset {\mathbb N}$, the set $A\cup F$ has not order 2 because
the set $\{ 2 10^{n} + 1 : n \in {\mathbb Z}^{+}\}$ has infinite many numbers.

We will now construct a subsequence of $\frac{(h-2)A(n)}{n} = \frac{A(n)}{n}$ in the following way:

\vspace{0.3cm}

$\{\frac{A(21)}{21}, \frac{A(201)}{201}, \frac{A(2001)}{2001}, \ldots, \frac{A(2 10^{n}+1)}{2 10^{n}+1}, \ldots \}$.

\vspace{0.3cm}

Let $A_n = \{ 2 \cdot 10^{n-1}+2, \ldots , 10^{n}\}$, $n \geq 2$. We then have

\begin{eqnarray*}
\displaystyle\sum_{i=2}^{n} |A_i | = \displaystyle\sum_{i=2}^{n} [10^{i}-(2 \cdot 10^{i-1}+2)]\\
=\displaystyle\sum_{i=2}^{n} (8 \cdot 10^{i-1}-2) = \displaystyle\sum_{i=1}^{n-1} (8 \cdot 10^{i}-2);
\end{eqnarray*}
hence,
\begin{eqnarray*}
\displaystyle\liminf_{n\rightarrow\infty} \frac{A(n)}{n}\leq
\displaystyle\lim_{n\rightarrow\infty}\frac{10 + \displaystyle\sum_{i=1}^{n-1} 
(8 \cdot 10^{i}-2)}{2 \cdot 10^{n}+1}=\frac{4}{9}.
\end{eqnarray*}

We next construct another subsequence of $\frac{A(n)}{n}$ given by:

\vspace{0.3cm}

$\{\frac{A(10)}{10}, \frac{A(100)}{100}, \ldots,  \frac{A(10^n)}{10^n}, \ldots\}$.

\vspace{0.3cm}

It is easy to see that
\begin{eqnarray*}
\displaystyle\limsup_{n\rightarrow\infty} \frac{A(n)}{n}\geq \displaystyle\lim_{n\rightarrow\infty}\frac{10 + 
\displaystyle\sum_{i=1}^{n-1} (8 \cdot 10^{i}-2)}{10^{n}}=\frac{8}{9}.
\end{eqnarray*}
Therefore, $\displaystyle\lim_{n\rightarrow\infty}\frac{A(n)}{n}$ does not exist. Then, $A$ is a 
finitely stable additive basis of order three
such that $\displaystyle\lim_{n\rightarrow\infty}\frac{(h-2)A(n)}{n}$ does not exist.

\section{Final Remarks}\label{sec3}

In this paper we have presented a constructive counterexample for Theorem 2.2 in
[Ferreira, L.A.: Finite Stable Additive Basis; Bull. Aust. Math. Soc.]. It seems
that the techniques proposed here can be used in the construction of additive basis of some order.

\begin{center}
\textbf{Acknowledgements}
\end{center}
This research has been partially supported by the Brazilian Agencies
CAPES and CNPq.

\end{document}